\documentclass[12 pt]{article}
\usepackage[letterpaper, margin = 1in]{geometry}

\usepackage{amsmath}
\usepackage{amsfonts}
\usepackage{amssymb}

\usepackage{graphicx}
\usepackage{epstopdf}
\usepackage{subcaption}

\usepackage[cal=boondox]{mathalfa}
\usepackage{hyperref}
\usepackage[noabbrev,capitalise]{cleveref}
\Crefname{ALC@unique}{Line}{Lines}
\crefformat{equation}{(#2#1#3)}

\author{Aurya Javeed\thanks{Center for Applied Mathematics, Cornell University, Ithaca, NY 14853.}}
\title{An Uncertainty Principle for Estimates of Floquet Multipliers}
\date{November 29, 2017}

\begin{document}
\maketitle

\begin{abstract}
  We derive a Cram\'er-Rao lower bound for the variance of Floquet multiplier estimates that have been constructed from stable limit cycles perturbed by noise. To do so, we consider perturbed periodic orbits in the plane. We use a periodic autoregressive process to model the intersections of these orbits with cross sections, then passing to the limit of a continuum of sections to obtain a bound that depends on the continuous flow restricted to the (nontrivial) Floquet mode. We compare our bound against the empirical variance of estimates constructed using several cross sections. The section-based estimates are close to being optimal. We posit that the utility of our bound persists in higher dimensions when computed along Floquet modes for real and distinct multipliers. Our bound elucidates some of the empirical observations noted in the literature; e.g., 
  \begin{enumerate}
    \item[(a)] it is the number of cycles (as opposed to the frequency of observations) that drives the variance of estimates to zero, and 
    \item[(b)] the estimator variance has a positive lower bound as the noise amplitude tends to zero.
  \end{enumerate}
\end{abstract}

\section{Introduction}

\begin{figure}\centering
  \includegraphics[width = 0.49\textwidth]{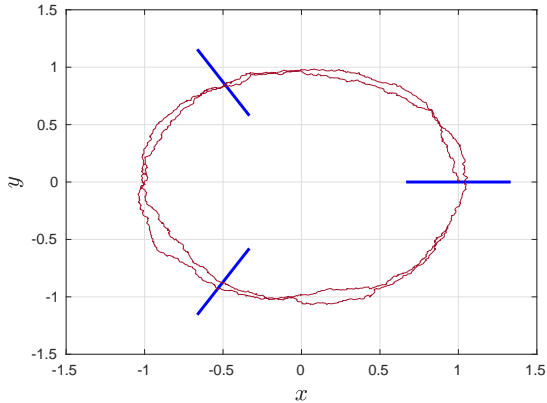}
  \caption{A stable periodic orbit perturbed by noise. Transverse cross sections are shown in blue. The amplitude of the noise is exaggerated for clarity.}\label{fig:c}
\end{figure}

This note studies systems like the one depicted in \cref{fig:c}. A stable periodic orbit of a dynamical system is perturbed by a small amount of noise so that the equations of motion become
\begin{align}\label{eq:em}
  \text{d}x = f(x)\text{d}t + g \text{d}w(t),
\end{align}
where $|g| \ll 1$ and $w(t)$ is an isotropic stochastic process. We derive an uncertainty principle governing the precision with which the Floquet multiplier of the unperturbed orbit can be specified from a noisy time series (e.g., the red one in \cref{fig:c}).

Recently, estimates of this nature have garnered attention in biomechanics \cite{g2014, holmes}. However, the question we answer in this note is fundamental:
\begin{quote}
  \textit{Given observations of a noisy dynamical process that is periodic when $g \to 0$, with how much precision can the stability of the $g = 0$ orbit be quantified?}
\end{quote}

In the language of statistical science, the uncertainty principle we derive is a \emph{Cram\'er-Rao lower bound}. Mathematically, these bounds are a consequence of the Cauchy-Schwarz inequality and read
\begin{align}\label{eq:cr}
  \text{var}_{\theta}(T) \geq J I(\theta)^{-1} J'.
\end{align}
Here, $\theta = (\theta_1, \ldots, \theta_n)$ is a set of parameters that characterize $\mathcal{p}_{\theta}$, the joint probability density from which the data are drawn. $T$ is a function of the data that estimates some value $v(\theta)$, $J$ is the Jacobian $[J]_i = \partial_{\theta_i}\mathbb{E}[T]$, and $I(\theta)$ is a matrix called the (Fisher) information. Under modest regularity conditions \cite{lc}\footnote{In \cite{lc}, \cref{eq:cr} is referred to as an information inequality (see Chapter 2, Section 5).}, $I(\theta)$ has entries
\begin{align*}
  [I(\theta)]_{ij} 
  = \mathbb{E} \left[ \Big( \partial_{\theta_i} \log \mathcal{p}_{\theta}(x) \Big)\Big(\partial_{\theta_j} \log \mathcal{p}_{\theta}(x) \Big) \right]
  = - \mathbb{E} \left[ \partial^2_{\theta_i\theta_j} \log \mathcal{p}_{\theta}(x) \right].
\end{align*}

The CR bound, \cref{eq:cr}, is significant because it holds for \emph{all} (sufficiently regular) estimators sharing a common expected value. 

\section{Derivation of the Uncertainty Principle}\label{sec:d}

To derive our uncertainty principle, we intersect the time series with $p$ transverse cross sections (see \cref{fig:c}). Ordering the intersections $x_i \in \mathbb{R}$ by time, we model $\{x_i\}_{i \geq 0}$ as a mean-zero periodic autoregressive process (PAR):
\begin{align}\label{eq:par}
  x_1 = \alpha_1x_0 + \epsilon_1, \qquad 
  x_2 = \alpha_2x_1 + \epsilon_2, \qquad 
  \cdots \qquad 
  x_p = \alpha_px_{p - 1} + \epsilon_p, \qquad
  \cdots
\end{align}
The indices of $\alpha$ are modulo $p$. For reasons that we will discuss later, no generality is lost in assuming \cref{eq:par} has mean zero. Also, we can assume that the random innovations $\epsilon_i$ are independent and identically distributed $N(0, \sigma^2)$ random variables. 

The derivation of our uncertainty principle consists of the following steps: First, we find the asymptotic variance of returns to each section under the PAR model, \cref{eq:par}. Then we use the variance to apply the CR bound to unbiased estimates of 
\begin{align*}
  \lambda = \prod_{k = 1}^p \alpha_k,
\end{align*}
the planar periodic orbit's Floquet multiplier. (This too is done assuming the PAR model.) Finally, we pass to the limit of a continuum of sections to obtain a lower bound that is independent of sections. Instead, the uncertainty principle we arrive at is a function of the linearized continuous flow of the noise-free system.

\subsection{Asymptotic Variance of Returns}\label{sec:av}

To derive the asymptotic variance of section returns, we suppose that the PAR begins at a known initial condition, $x_0$, on Section $0$. The process then evolves forward as 
\begin{align*}
\begin{split}
  x_1 &= \alpha_1x_0 + \epsilon_1\\
  x_2 &= \alpha_2x_1 + \epsilon_2 = \alpha_2\alpha_1x_0 + \epsilon_2 + \alpha_2\epsilon_1\\
  x_3 &= \alpha_3x_2 + \epsilon_3 = \alpha_3\alpha_2\alpha_1x_0 + \epsilon_3 + \alpha_3\epsilon_2 + \alpha_3\alpha_2\epsilon_1 \\
&\enspace\vdots
\end{split}
\end{align*}
After $p$ steps, the PAR returns to the starting section at the point
\begin{align*}
  x_p = \lambda x_0 + \epsilon_p + \sum_{i = 1}^{p - 1} \epsilon_i \prod_{j = i}^{p - 1} \alpha_{j + 1}
\end{align*}
(where the summation is empty if $p = 1$). Reindexing $j$ as $j - 1$ then implies
\begin{align*}
  \text{var}(x_p) = \sigma^2 \left( 1 + \sum_{i = 2}^p \prod_{j = i}^p \alpha_j^2 \right),
\end{align*}
and it follows that
\begin{align*}
  \text{var}(x_{p + np}) = \text{var}(x_p)\sum_{\ell = 0}^n (\lambda^2)^{\ell}.
\end{align*}
In the $n\to\infty$ limit, this variance becomes $\text{var}(x_p)(1 - \lambda^2)^{-1}$.

If the process instead started on Section $k$, the same reasoning implies that $\text{var}(x_{k + p})$ will be $\text{var}(x_p)$, except with the indices of $\alpha$ shifted from $j$ to $j + k$. As $n \to \infty$, the section on which the process started is essentially insignificant, hence the asymptotic variance of returns to Section $k$ is
\begin{align}\label{eq:av}
  \lim_{n \to \infty} \text{var}(x_{(k + p) + np})
  = \text{var}(x_{k + p}) \left[ \lim_{n \to \infty} \sum_{\ell = 0}^n (\lambda^2)^{\ell} \right]
  = \frac{\sigma^2}{1 - \lambda^2} \left( 1 + \sum_{i = 2}^p \prod_{j = i}^p \alpha_{j+k}^2 \right ).
\end{align}

\subsection{Applying the CR Bound}

Having obtained the asymptotic variance of returns, we now apply the CR bound to unbiased estimators of $\lambda$; i.e., to estimators having $\mathbb{E}[T] = \lambda = v(\theta)$ for $\theta = (\alpha_1, \ldots, \alpha_p)$. 

The Jacobian row vector in \cref{eq:cr} has entries $[J]_i = \lambda/\alpha_i$. To derive the information, we note that
\begin{align*}
  \mathcal{p}_{\theta}(x_k| x_{k - 1}) = \frac{1}{\sqrt{2\pi\sigma^2}} \exp\left[ -\frac{(x_k - \alpha_kx_{k - 1})^2}{2\sigma^2} \right],
\end{align*}
yielding
\begin{align*}
  \mathcal{p}_{\theta}(x_0, x_1, \ldots, x_p) = \prod_{k = 1}^p \mathcal{p}_{\theta}(x_k | x_{k - 1}) 
= \prod_{k = 1}^p \frac{1}{\sqrt{2\pi\sigma^2}} \exp\left[ -\frac{(x_k - \alpha_kx_{k - 1})^2}{2\sigma^2} \right].
\end{align*}
Thus, for $p$ successive observations, $I(\theta)$ is asymptotically diagonal with entries
\begin{align*}
  [I(\theta)]_{kk} = - \mathbb{E}\left[ \partial^2_{\alpha_k \alpha_k}\log \mathcal{p}_{\theta} \right] = \frac{\text{var}\left( x_{k - 1} \right)}{\sigma^2} = \frac{1}{1 - \lambda^2} \left( 1 + \sum_{i = 2}^p \prod_{j = i}^p \alpha_{j + k - 1}^2 \right ),
\end{align*}
where we have used \cref{eq:av}. Consequently,
\begin{align*}
  JI(\theta)^{-1}J' 
  = \sum_{k = 1}^p\frac{\lambda^2(1 - \lambda^2)}{\left( \alpha_k^2 + \alpha_k^2\sum_{i = 2}^p \prod_{j = i}^p \alpha_{j + k - 1}^2 \right )} 
  = \lambda^2(1 - \lambda^2)\sum_{k = 1}^p \left(\sum_{i = 1}^p \prod_{j = i}^p \alpha_{j + k}^2 \right )^{-1}.
\end{align*}
This last equality follows from reindexing (cf. the reindex in the first paragraph of \cref{sec:av}).

\subsection{Continuum Limit}

Thus far we have assumed nothing about the distance between sections, so we are justified in spacing them so that the same amount of time, $\Delta t$, elapses between $x_i$ and $x_{i + 1}$, regardless of the value of $i$. By writing
\begin{align*}
  JI(\theta)^{-1}J'
  = \lambda^2(1 - \lambda^2)\sum_{k = 1}^p \Delta t \left( \sum_{i = 1}^p \Delta t \prod_{j = i}^p \alpha_{j + k}^2 \right )^{-1},
\end{align*}
we find ourselves in a position to take the continuum limit $(p, \Delta t) \to (\infty, 0)$ with $p\Delta t$ held constant. The sums will clearly limit to integrals, but we devote a few sentences to what becomes of $\prod_{j = i}^p \alpha^2_{j + k}$, which is a mapping from Section $k + i$ forward to $k + p \equiv k$. 

Suppose $\gamma$ denotes the unperturbed periodic orbit, and let $\gamma(t = 0)$ be the intersection of Section 0 and $\gamma$. For the $f$ in \cref{eq:em}, we define $\Phi(t)$ to be the principal solution matrix of
\begin{align*}
  \frac{d\xi}{dt} = Df|_{\gamma(t)}\xi,
\end{align*}
which is known as the variational equation of $\gamma$. The matrix $\Phi$ has a set of invariant subspaces that are independent of $\gamma(0)$ and are called the Floquet modes of $\gamma$ \cite{chicone}. One of the modes is always trivial, traced out by the tangent vector $f \circ \gamma(t)$. There is also a nontrivial mode corresponding to $\lambda$. We let $\phi:[0, \infty) \to \mathbb{R}$ be the restriction of $\Phi$ to this nontrivial mode. If the sections are oriented appropriately, the transition map from Section $k + i$ to Section $k + p$ follows as
\begin{align*}
  \phi(b + \tau)\phi(b + a)^{-1},
\end{align*}
where $\tau$ is the period of $\gamma$, and quantities $a$ and $b$ are the times that $\gamma$ first intersects Sections $i$ and $k$, respectively.

Returning to the limit of a continuum of sections, we therefore have that
\begin{align*}
  JI(\theta)^{-1}J'
  = \lambda^2(1 - \lambda^2)\int_0^{\tau} ds \left(\int_0^{\tau} dt \left[\phi(s + \tau)\phi(t + s)^{-1} \right]^2 \right )^{-1},
\end{align*}
where $(a, b)$ has been relabeled as $(t, s)$. Floquet's theorem \cite{chicone} simplifies this equality to
\begin{align*}
  JI(\theta)^{-1}J' = (1 - \lambda^2)\int_0^{\tau} ds \left(\int_0^{\tau} dt \left[\phi(s)\phi(t + s)^{-1} \right]^2 \right )^{-1}.
\end{align*}
Finally, we recall that the information in this expression corresponds to one cycle of the perturbed orbit. For $n$ cycles, $I(\theta)$ will be $n$ times as large, bringing us to our main result:
\begin{align}\label{eq:up}
  \boxed{\text{var}(T) \geq \frac{(1 - \lambda^2)}{n}\int_0^{\tau} \frac{ds}{\phi(s)^2\int_0^{\tau} \phi(t + s)^{-2} dt}.}
\end{align}

Though this uncertainty principle was derived by considering a planar periodic orbit, we expect \cref{eq:up} to persist for real and distinct multipliers of higher-dimensional periodic orbits when $\phi$ is the restriction of $\Phi$ to the Floquet mode of the multiplier in question. Why do we expect this? Because the best one can hope for when estimating a multiplier is that the entire time series resides in the invariant subspace associated with the multiplier.

In the next section, we show that our uncertainty principle and its generalization to higher-dimensional orbits are indeed reasonable. However, before doing so, we fulfill our promise to explain why no generality was lost when assuming the PAR has mean zero and that $\epsilon_i \overset{i.i.d.}{\sim} N(0, \sigma^2)$. 
\begin{itemize}
  \item For the mean-zero assumption, note that if \cref{eq:par} had nonzero mean, it could be centered by the transformation $x_i \mapsto (x_i - \mathbb{E}[x_i])$, with $\mu_{i\, \text{mod}\, p} = \mathbb{E}[x_i]$ defining an additional parameter to estimate. But this parameter does not influence our uncertainty principle since $\partial_{\mu_i} \lambda = 0$.

  \item Regarding the $\epsilon_i$ assumption: if the stochastic process driving the noisy equations of motion is a Wiener process, we expect the innovations near the continuum limit to be i.i.d. $N(0, \sigma^2)$ (cf. the Euler-Maruyama algorithm, which converges as $\Delta t \to 0$). Generalizations to other noise types are briefly discussed in \cref{sec:c}.
\end{itemize}

\section{Numerical Validation}

\begin{figure}\centering
  \begin{subfigure}{.5\textwidth}\centering
    \includegraphics[width=0.9\linewidth]{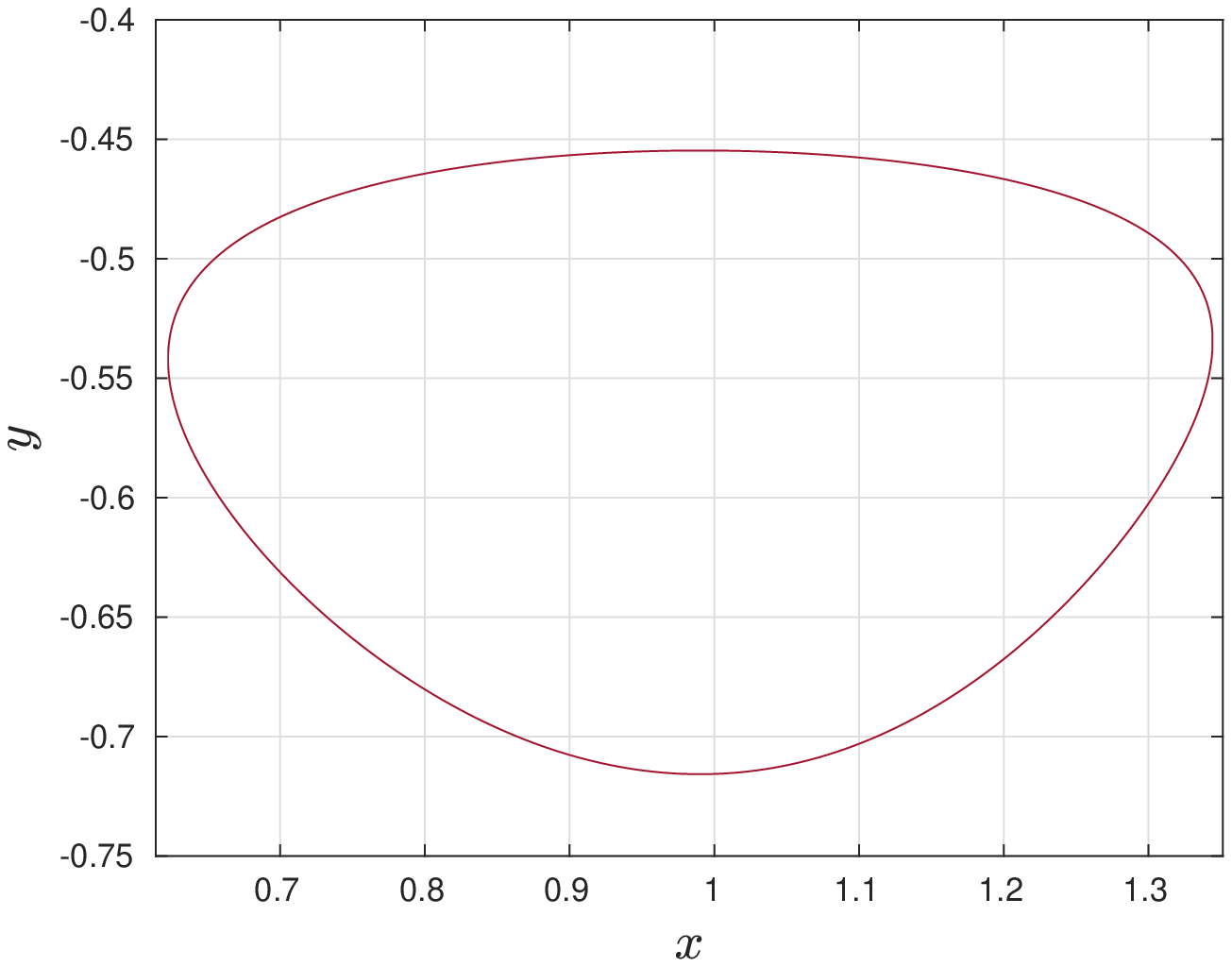}\caption{}\label{fig:vdp}
  \end{subfigure}%
  \begin{subfigure}{.5\textwidth}\centering
    \includegraphics[width=0.9\linewidth]{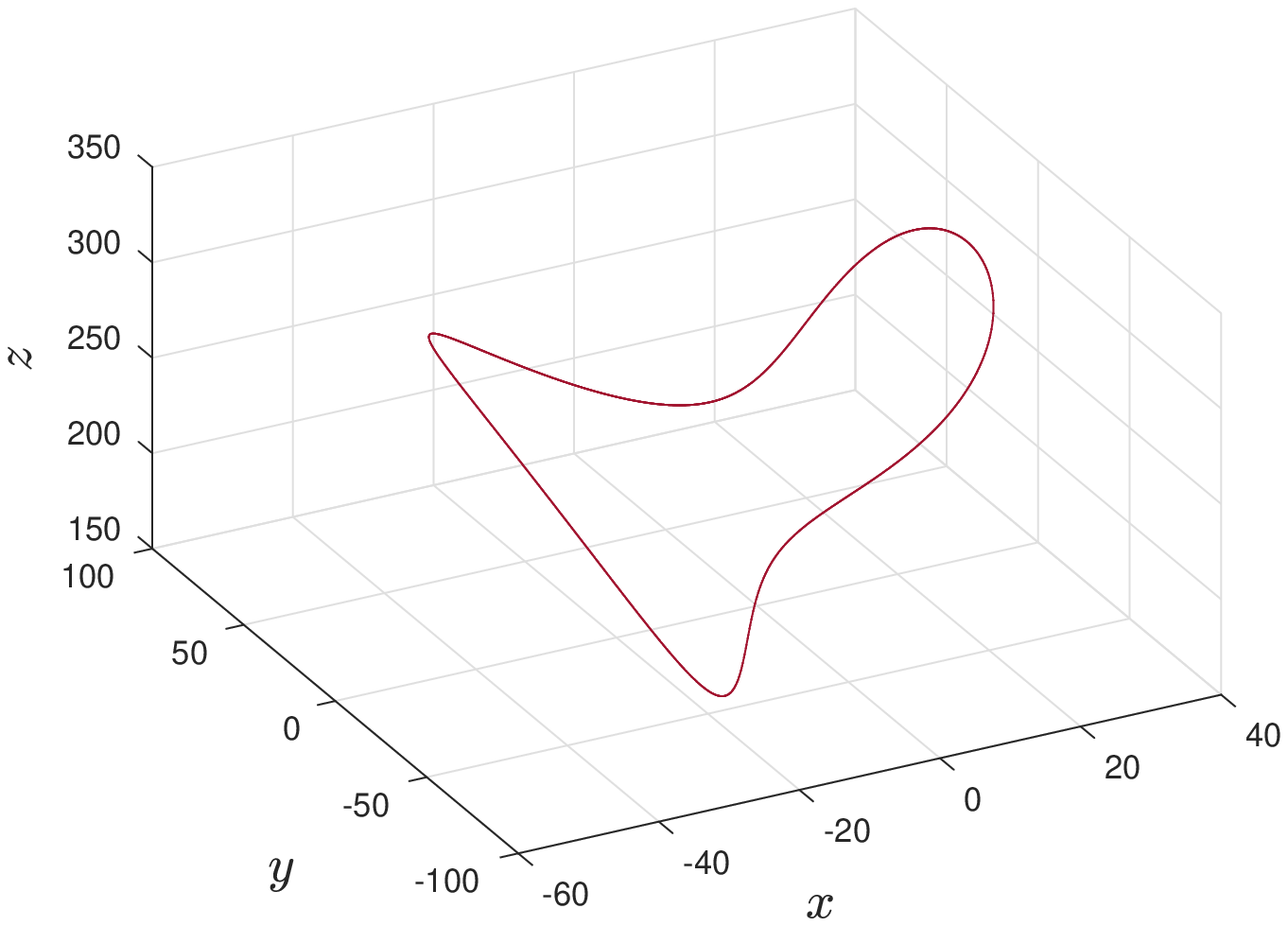}\caption{}\label{fig:l}
  \end{subfigure}
  \caption{Periodic orbits of (a) the van der Pol system, \cref{eq:vdp}, and (b) the Lorenz system, \cref{eq:l}.}\label{fig:o}
\end{figure}

We consider periodic orbits from two different dynamical systems. Both of them are shown in \cref{fig:o}. 

The first orbit (\cref{fig:vdp}) belongs to the van der Pol system
\begin{align}\label{eq:vdp}
\begin{split}
  \varepsilon\dot{x} &= y - \frac{x^3}{3} + x \\
  \dot{y} &= a - x,
\end{split}
\end{align}
with $(\varepsilon, a) = (0.1, 0.99)$. In the leftmost columns of \cref{tbl:r}, we record the multiplier of this orbit, along with the square root of our uncertainty principle when $n = 100$. Next to these two columns are statistics of 100 realizations of a numerical simulation. In each realization, the equations of motion are perturbed by Gaussian noise having amplitude $5\times 10^{-5}$, and the multiplier of the orbit is estimated from 100 sample path returns to 50 sections. The estimation is performed using least squares, fitting scalars $\alpha_i$ from section to section.

The predictive power of our uncertainty principle is striking! Additionally, it suggests that the section-based estimator is quite reasonable.

\begin{table}
  \caption{Numerical validation of our uncertainty principle, \cref{eq:up}. For each of the orbits in \cref{fig:o}, the multiplier and square root of our UP are tabulated alongside statistics of 100 independent numerical simulations that estimate $\lambda$ from the perturbed system's returns to several sections.}\label{tbl:r}
  \centering\small
  \begin{subtable}{.5\linewidth}\centering
    {\begin{tabular}[t]{ c c | c c}
      $\lambda$ & $\sqrt{\text{UP}}$ & \verb|mean(| $\hat{\lambda}$ \verb|)| & \verb|std(| $\hat{\lambda}$ \verb|)|\\
	  \hline
	  0.3854 & 0.0532 & 0.3953 & 0.0541\\
    \phantom{a}
    \end{tabular}}\caption{Results for the van der Pol orbit.}
  \end{subtable}%
  \begin{subtable}{.5\linewidth}\centering
    {\begin{tabular}[t]{c | c c | c c}
      & $\lambda_i$ & $\sqrt{\text{UP}}$ & \verb|mean(| $\hat{\lambda}_i$ \verb|)| & \verb|std(| $\hat{\lambda}_i$ \verb|)|\\
	  \hline
	  $i = 1$ & -0.6162 & 0.0606 & -0.6157 & 0.0749 \\
	  2    & -0.0026 & 0.0009 & -0.0031 & 0.0014
    \end{tabular}}\caption{Results for the Lorenz orbit.}
  \end{subtable}
\end{table}

The right half of \cref{tbl:r} presents these same quantities for the second periodic orbit we consider (\cref{fig:l}). This second orbit is one of two present the Lorenz system
\begin{align}\label{eq:l}
\begin{split}
  \dot{x} &= \sigma (y - x) \\
  \dot{y} &= rx - y - xz    \\
  \dot{z} &= -bz + xy,
\end{split}
\end{align} 
when parameters $(\sigma,\, r,\, b) = (10,\, 240,\, 8/3)$. Per \cref{sec:d}, the right-hand side of our uncertainty principle was computed one multiplier at a time. In simulations, the amplitude of the noise perturbing the orbit  was increased to $3 \times 10^{-2}$, following \cite{g2014}. Even in this higher dimensional setting, the agreement between our uncertainty principle and the empirical results is very encouraging.\footnote{Interestingly, the eigenvalues of our monodromy matrix are conditioned oppositely: $\lambda_1$ is more robust to perturbations of the matrix (see Section 7.2.2 of \cite{gvl}); but---as anticipated by our uncertainty principle---estimates of this slow multiplier are more variable.}

\section{Concluding Remarks}\label{sec:c}

This note presents a fundamental and broadly-applicable result:
\begin{quote}
  \textit{Given a time series of a perturbed vector field with a periodic orbit, we have characterized an intrinsic precision with which the stability of the deterministic periodic orbit can be determined.}
\end{quote}
Insofar as generality is concerned,
\begin{itemize}
  \item CR bounds also apply to biased estimators. If the bias is small (as it is in \cref{tbl:r}), the adjustment to our uncertainty principle will be small too \cite{lc}.
  
  \item We expect our uncertainty principle to persist in higher dimensions when the multipliers of the periodic orbit are real and distinct. To evaluate the bound in this setting, we take $\phi$ to be the restriction of $\Phi$ to the Floquet mode corresponding to the multiplier in question. Our numerical results in \cref{tbl:r} suggest that bounds obtained this way may not be unreasonably loose. 
  
  \item It should be possible to extend our uncertainty principle to other types of noise by adjusting the information accordingly (e.g., if instead of Brownian motion, the perturbed equations of motion are driven by a process with isotropic Laplace increments, $I(\theta)$ will be twice as large).
\end{itemize}

Finally, we remark that our uncertainty principle is significant because of the quantities that do and do not appear in it. For example, there is no dependence on the amplitude, $g$, of the perturbing noise (though, of course, it must be small enough for a linear approximation of the transverse dynamics to be valid). In addition, our uncertainty principle depends on the number of cycles but not on the frequency of observations, meaning the former and not the latter drives the variability of estimates to zero.
\begin{quote}
  \textit{This latter observation is particularly significant because it implies that if a spectator (e.g., a control) is trying to infer multipliers from a noisy periodic orbit in ``stochastic equilibrium" \cite{cyclo}, they will be unable to do so to an arbitrary degree of precision in a fixed amount of time.}
\end{quote}

The code used in this note is happily provided upon request. It uses \cite{viswanath} and Chebfun \cite{chebfun} to minimize numerical error.

\section*{Acknowledgements}

A sincere thanks to Professors John Guckenheimer and Giles Hooker. This work would not have been possible without their guidance.

\bibliographystyle{plain}
\bibliography{bib}

\end{document}